\documentclass[12pt]{amsart} 
\usepackage{amssymb,amsmath,latexsym, enumerate} 

\usepackage{color}

\numberwithin{equation}{section}
\theoremstyle{plain} 
\newtheorem{thm}{Theorem}[section]

\newtheorem{lem}[thm]{Lemma}
\newtheorem{rmk}[thm]{Remark}
\newtheorem{ex}[thm]{Example}
\newtheorem{cor}[thm]{Corollary}

\newtheorem{ques}[thm]{Question}

\theoremstyle{definition}

\begin{document}
\newcommand{\hgt}{\operatorname{ht}}
\newcommand{\ext}{\operatorname{Ext}}
\newcommand{\Hom}{\operatorname{Hom}}
\newcommand{\extx}{\operatorname{Ext}_{R_x}}

\newcommand{\exti}{\operatorname{Ext}^i_R}
\newcommand{\extijM}{\operatorname{Ext}^i_R(R/I,H^j_I(M))}
\newcommand{\extjiM}{\operatorname{Ext}^j_R(R/I,H^i_I(M))}

\newcommand{\extijiM}{\operatorname{Ext}^i_R(R/I,H^{j+1}_I(M))}
\newcommand{\extijR}{\operatorname{Ext}^i_R(R/I,H^j_I(R))}
\newcommand{\extidR}{\operatorname{Ext}^i_R(R/I,H^d_I(R))}
\newcommand{\extidiR}{\operatorname{Ext}^i_R(R/I,H^{d-1}_I(R))}

\newcommand{\hjR}{H^j_I(R)}
\newcommand{\hiR}{H^i_I(R)}
\newcommand{\hjM}{H^j_I(M)}
\newcommand{\hjiM}{H^{j+1}_I(M)}
\newcommand{\hiM}{H^i_I(M)}
\newcommand{\hhM}{H^h_I(M)}
\newcommand{\hnM}{H^n_I(M)}
\newcommand{\hniM}{H^{n-1}_I(M)}

\newcommand{\hjRx}{H^j_I{R_x}(R_x)}

\newcommand{\hdM}{H^d_I(M)}
\newcommand{\hdR}{H^d_I(R)}
\newcommand{\hdiM}{H^{d-1}_I(M)}
\newcommand{\hdiR}{H^{d-1}_I(R)}
\newcommand{\hdiiR}{H^{d-2}_I(R)}
\newcommand{\hiIM}{H^i_I(M)}
\newcommand{\hgIM}{H^g_I(M)}
\newcommand{\supp}{\operatorname{Supp}_R}
\newcommand{\ass}{\operatorname{Ass}_R}
\newcommand{\suppx}{\operatorname{Supp}_{R_x}}
\newcommand{\assx}{\operatorname{Ass}_{R_x}}
\newcommand{\depth}{\operatorname{depth}}
\newcommand{\ch}{\operatorname{v}}
\newcommand{\ann}{\operatorname{ann}}
\newcommand{\soc}{\operatorname{soc}}
\newcommand{\starsoc}{^*\operatorname{soc}}

\title[Local cohomology modules]{Local
cohomology modules with infinite dimensional socles}
\author{Thomas Marley}
\address{University of Nebraska-Lincoln\\
Department of Mathematics and Statistics\\
Lincoln, NE 68588-0323}
\email{tmarley@math.unl.edu}
\urladdr{http://www.math.unl.edu/\textasciitilde tmarley}
\author{Janet C. Vassilev}
\address{University of Arkansas\\
Department of Mathematical Sciences\\
Fayetteville, AR 72701}
\email{jvassil@uark.edu}
\urladdr{http://comp.uark.edu/\textasciitilde jvassil}

\thanks{The first author was partially supported by NSF grant DMS-0071008.}
\date{August 27, 2002}

\begin{abstract}  In this paper we prove the following generalization of a result of Hartshorne:
Let $T$ be a commutative Noetherian local ring of dimension
at least two, $R=T[x_1,\dots,x_n]$, and $I=(x_1,\ldots,x_n)$.  Let $f$ be a homogeneous element of $R$ such that the coefficients of $f$ form a system of parameters for $T$.  Then the socle of
$H^n_I(R/fR)$ is infinite dimensional.
\end{abstract}

\maketitle

\section{Introduction}

The third of Huneke's four problems in
local cohomology \cite{Hu} is to determine when $H^i_I(M)$ is Artinian for a given ideal $I$ of a commutative Noetherian local ring $R$ and finitely generated
$R$-module $M$.   An $R$-module $N$ is Artinian if and only $\supp N\subseteq \{m\}$ and 
$\Hom_R(R/m,N)$ is finitely generated, where $m$ is the maximal ideal of $R$.  Thus, Huneke's 
problem may be separated into two subproblems: 

\begin{itemize}
\item When is $\supp H^i_I(M) \subseteq \{m\}$?
\item When is $\Hom_R(R/m,H^i_I(M))$ finitely generated?
\end{itemize}

This article is concerned with the second question.
For an $R$-module $N$, one may identify $\Hom_R(R/m,N)$ with the submodule $\{x\in N\mid mx=0\}$,
which is an $R/m$-vector space
called the {\it socle} of $N$ (denoted $\soc_RN$).
It is known that if $R$ is an unramified regular local ring then the local cohomology modules $H^i_I(R)$ have finite dimensional socles for all $i\ge 0$ and all ideals $I$ of $R$ (\cite{HS}, \cite{L1}, \cite{L2}).   The first example of a local cohomology module  with an infinite dimensional socle was given in 1970 by Hartshorne \cite{Ha}:
Let $k$ be a field, $R=k[[u,v]][x,y]$, $P=(u,v,x,y)R$, $I=(x,y)R$, and $f=ux+vy$.  Then $\soc_{R_P} H^2_{IR_P}(R_P/fR_P)$ is infinite dimensional.  Of course, since $I$ and $f$ are homogeneous, this is equivalent to saying
that $\Hom_R(R/P,H^2_I(R/fR))$ (the {\it $^*$socle} of $H^2_I(R/fR)$) is infinite dimensional.  Hartshorne proved this by exhibiting an infinite set of linearly independent elements in
the $^*$socle of $H^2_I(R)$.

In the last 30 years there have been few results in the literature which explain or generalize
 Harthshorne's example.   For affine semigroup rings, a remarkable result proved
by Helm and Miller \cite{HM} gives necessary and sufficient
conditions (on the semigroup) for the ring to possess a local cohomology module (of a finitely
generated module) having infinite dimensional socle.   Beyond that work, however, little has
been done.

In this paper we prove the following:

\begin{thm} \label{mainthm} Let $(T,m)$ be a Noetherian local of dimension at least two.  Let $R=T[x_1,\dots,x_n]$
be a polynomial ring in $n$ variables over $T$, $I=(x_1,\dots,x_n)$, and $f\in R$ a homogeneous
polynomial whose coefficients form a system of parameters for $T$.  Then the $^*$socle of
$H^n_I(R/fR)$ is infinite dimensional.
\end{thm}

Hartshorne's example is obtained by letting $T=k[[u,v]]$, $n=2$, and $f=ux+vy$ (homogeneous of
degree 1).  Note, however, that we do not require the coefficient ring to be regular, or even
Cohen-Macaulay.  As a further illustration, consider the following:

\begin{ex} {\rm Let $R=k[[u^4,u^3v,uv^3,v^4]][x,y,z]$, $I=(x,y,z)R$, and $f=u^4x^2+v^8yz$.  Then the
$^*$socle of
$H^3_I(R/fR)$ is infinite dimensional.} 
\end{ex}

Part of the proof of Theorem \ref{mainthm} was inspired by the recent work of Katzman \cite{Ka}
where information on 
the graded pieces of $H^n_I(R/fR)$ is obtained by examining matrices
of a particular form.  We apply this technique in the proof of Lemma \ref{keylemma}.

Throughout all rings are assumed to be commutative with identity.
The reader should consult \cite{Mat} or \cite{BH} for any unexplained terms or notation and
\cite{BS} 
for the basic properties of local cohomology.

\section{The Main Result}

Let $R=\oplus R_{\ell}$ 
be a Noetherian ring graded by the nonnegative integers.   Assume
$R_0$ is local and let $P$ be the homogeneous maximal ideal of $R$. 
Given a finitely
generated graded
$R$-module $M$ we define the \it $^*$socle  \rm of $M$ by
\begin{align}
\starsoc_RM &= \{x\in M \mid Px=0\} \notag \\
&\cong \Hom_R(R/P,M). \notag
\end{align}

Clearly, $\starsoc_R M\cong \soc_{R_P}M_P$.  An interesting special case of Huneke's third problem
is the following:

\begin{ques} \label{Q1} Let $n:=\mu_R(R_+/PR_+)$, the minimal number of generators of $R_+$. When
 is $\starsoc H^n_{R_+}(R)$ finitely generated?
\end{ques}

For $i\in \mathbb N$ it is well known that $H^i_{R+}(R)$ is a graded $R$-
module, each graded piece $H^i_{R_+}(R)_{\ell}$
is a finitely generated $R_0$-module, and $H^i_{R_+}(R)_{\ell}=0$ for all sufficiently large
integers $\ell$ (\cite[15.1.5]{BS}).  
If we know \it a priori \rm that $H^n_{R+}(R)_{\ell}$ has finite length for all
$\ell$ (e.g., if $\supp H^n_{R_+}(R)\subseteq \{P\}$), then Question \ref{Q1} is equivalent to:

\begin{ques} \label{Q2} When is $\Hom_R(R/R_+,H^n_{R_+}(R))$ finitely 
generated?  
\end{ques}

We give a partial answer to these questions for hypersurfaces.
For the remainder of this section we
adopt the following notation:  Let $(T,m)$ be a local ring of dimension $d$ and 
$R=T[x_1,\dots,x_n]$ a polynomial ring in $n$ variables over $T$.  We endow $R$ with an
$\mathbb N$-grading by setting $\deg T=0$ and $\deg x_i=1$ for all $i$.  Let $I=R_+=(x_1,\dots,x_n)R$ and $P=m+I$ the homogeneous maximal ideal of $R$.  Let $f\in R$ be a homogeneous element of
degree $p$ and $C_f$ the ideal of $T$ generated by the nonzero coefficients of $f$.

Our main result is the following:

\begin{thm} \label{mainresult} Assume $d\ge 2$ and the (nonzero) coefficients of $f$ form a
system of parameters for $T$.  Then $\starsoc_R H^n_I(R/fR)$ is not finitely generated.
\end{thm}

The proof of this theorem will be given in a series of lemmas below.  Before proceeding with the
proof we make a couple of remarks:

\begin{rmk}{\rm

\begin{enumerate}[(a)]
\item If $d\le 1$ in Theorem \ref{mainresult} then $\starsoc H^n_I(R/fR)$ is 
finitely generated.  This follows from \cite[Corollary 2]{DM} since $\dim R/I=\dim T\le 1$.
\item The hypothesis that the nonzero coefficients of $f$ form a system of parameters for $T$ is
stronger than our proof requires.  One only needs that
$C_f$ be $m$-primary and that there exists a dimension 2 ideal
containing all but two of the coefficients of $f$. (See the proof of Lemma \ref{keylemma}.)
\end{enumerate}}
\end{rmk}

The following lemma identifies the support of $H^n_I(R/fR)$ for a homogeneous element $f\in R$.   
This lemma also follows from a much more general result recently proved by Katzman
and Sharp \cite[Theorem 1.5]{KS}.

\begin{lem} \label{support} Let $f\in R$ be a homogeneous element.  Then
$$\operatorname{Supp}_RH^n_I(R/fR)=\{Q\in \operatorname{Spec}R \mid Q\supseteq I+C_f\}.$$
\end{lem}

{\it Proof:}  It is enough to prove that $H^n_I(R/fR)=0$ if and only if $C_f=T$.
As $H^n_I(R/fR)_k$ is a finitely generated $T$-module
for all $k$, we have by Nakayama that  $H^n_I(R/fR)=0$ if and only if $H^n_I(R/fR)\otimes_T
T/m=0$.  Now
\begin{align}
H^n_I(R/fR)\otimes_T T/m&\cong H^n_I(R/fR\otimes_T T/m) \notag \\
&\cong H^n_{N}(S/fS) \notag
\end{align}
where $S=(T/m)[x_1,\dots,x_n]$ is a polynomial
ring in $n$ variables over a field and $N=(x_1,\dots,x_n)S$.
As $\dim S=n$, we see that $H^n_N(S/fS)=0$
if and only if the image of $f$
modulo $m$ is nonzero.  Hence, $H^n_I(R/fR)=0$ if and only if at least one coefficient of 
$f$ is a unit, i.e., $C_f=T$.
\qed

We are mainly interested in the case the coefficients of $f$ generate an $m$-primary ideal:

\begin{cor} \label{supportmax} Let $f\in R$ be homogeneous and suppose $C_f$ is $m$-primary.  Then
$$\operatorname{Supp}_R H^n_I(R/fR)=\{P\}.$$ 
\end{cor}

Our next lemma is the key technical result in the proof of Theorem \ref{mainresult}.

\begin{lem} \label{matrixlemma}  Suppose $u,v\in T$ such that $\operatorname{ht}
(u,v)T=2$.  For each integer $n\ge 1$ let $M_n$ be the cokernel of 
$\phi_n:T^{n+1}\rightarrow T^n$ where $\phi_n$ is represented by the matrix
$$
A_n=\begin{pmatrix}
u & v & 0 &0 & \cdots & 0 & 0 \\
0 & u & v & 0& \cdots & 0 & 0 \\
0 & 0 & u & v & \cdots & 0 & 0  \\
\vdots & \vdots & \vdots  &\ddots & \ddots & \vdots & \vdots \\
0 & 0 & 0 & 0 & \cdots & u & v \\
\end{pmatrix}_{n\times (n+1)}.
$$
Let $J=\cap_{n\ge 1}\ann_TM_n$.  Then $\dim T/J=\dim T$.
\end{lem}

{\it Proof:} Let $\hat T$ denote the $m$-adic completion of $T$.
Then $\operatorname{ht} (u,v)\hat T=2$, $\ann_TM_n = \ann_{\hat T}(M_n\otimes_T\hat{T}) \cap T$,
and $\dim T/(I\cap T)\ge \dim \hat T/I$ for all ideals $I$ of
$\hat T$.  Thus,
we may assume $T$ is complete.  Now let $p$ be a prime ideal of $T$ such that
$\dim T/p=\dim T$.  Since $T$ is catenary, $\operatorname{ht}(u,v)T/p=2$.  
Assume the lemma is true for complete domains.  Then 
$\cap_{n \ge 1}\ann_{T/p}(M_n\otimes_TT/p)=p/p$.
Hence
\begin{align}
J&=\cap_{n\ge 1}\ann_TM_n \notag \\
&\subseteq \cap_{n\ge 1}\ann_T (M_n\otimes_T T/p) \notag \\
&= p, \notag
\end{align}
which implies that $\dim T/J\ge \dim T/p=\dim T$.
Thus, it suffices to prove the lemma for complete domains. 

As $T$ is complete, the integral closure $S$ of $T$ is a finite $R$-module.  Since 
$\operatorname{ht} (u,v)S=2$
(\cite[Theorem 15.6]{Mat}) and
$S$ is normal, $\{u,v\}$ is a regular sequence on $S$.  It is easily seen that $I_n(A_n)$, the
ideal of $n\times n$ minors of $A_n$, is $(u,v)^nT$.  
By the main result of \cite{BE} we obtain $\ann_S (M_n\otimes_TS)=
(u,v)^nS$.  Hence $\ann_TM_n\subseteq (u,v)^nS\cap T$.  As $S$ is a finite $T$-module there exists
an integer $k$ such that $\ann_TM_n\subseteq (u,v)^{n-k}T$ for all $n\ge k$.  Therefore,
$\cap_{n\ge 1}\ann_TM_n=(0)$, which completes the proof.
\qed

\begin{lem} \label{keylemma} Assume $d\ge 2$ and let $f\in R$ be a homogeneous element of
degree $p$ such that 
the coefficients of $f$ form a system of parameters for $T$.  Then 
 $\dim T/\ann_T H^n_I(R/fR) \ge 2$.
\end{lem}

{\it Proof:} Let $c_1,\dots,c_d$ be the nonzero coefficients of $f$.  Let $T'=T/(c_3,\dots,c_d)T$
and $R'=T'[x_1,\dots,x_n]\cong R/(c_3,\dots,c_d)R\cong R\otimes_TT'$.  Since 
\begin{align}
\dim T/\ann_T H^n_I(R/fR) &\ge
\dim T/ \ann_T (H^n_I(R/fR)\otimes_TT') \notag \\
&= \dim T'/\ann_{T'} H^n_{IR'}(R'/fR'), \notag
\end{align}
we may assume that $\dim T=2$ and $f$ has exactly two nonzero terms.

For any  $w\in R$ there is a surjective map 
$H^n_I(R/wfR)\to H^n_I(R/fR)$.  Hence, $\ann_T H^n_I(R/wfR) \subseteq \ann_T H^n_I(R/fR)$.
Thus, we may assume that the  terms of $f$ have no (nonunit) common factor.  Without loss
of generality, we may write $R=T[x_1,\dots,x_k,y_1,\dots y_r]$ and
$f=ux_1^{d_1}\cdots x_k^{d_k} + vy_1^{e_1}\cdots y_r^{e_r}=u\mathbf x^{\mathbf d}+v\mathbf y
^{\mathbf e}$, where
$\{u,v\}$ is a system of parameters for $T$.  As $f$ is homogeneous,
 $p=\sum_id_i=\sum_i e_i$.

Applying the right exact functor $H^n_I(\cdot)$ to 
 $R(-p)\xrightarrow{f} R \to R/fR \to 0$ we obtain the
exact sequence
$$H^n_I(R)_{-\ell-p} 
\xrightarrow{f} H^n_I(R)_{-\ell} \to H^n_I(R/fR)_{-\ell}\to 0$$
for each $\ell\in \mathbb Z$.  For each $\ell$, $H^n_I(R)_{-\ell}$ is a free $T$-module
with basis 
$$\{\mathbf x^{-\mathbf \alpha}\mathbf y^{-\mathbf \beta} \mid \sum_{i,j}\alpha_i+\beta_j
=\ell, \alpha_i>0, \beta_j>0 \ \forall \ i,j\}$$
(e.g., \cite[Example 12.4.1]{BS}).
Let $q$ be an arbitrary positive integer and let $\ell(q)=qp+k+r$.  Define $L_{-\ell(q)}$ 
to be the free $T$-summand
of
$H^n_I(R)_{-\ell(q)}$ spanned by the set 
$$\{\mathbf x^{-s\mathbf d-\mathbf 1}\mathbf y^{-t\mathbf e - \mathbf 1}
\mid s+t=q, s,t\ge 0\}.$$
Then the cokernel of $\delta_q:L_{-\ell(q+1)}\xrightarrow{f} L_{-\ell(q)}$ is a direct summand 
(as a $T$-module) of 
$H^n_I(R/fR)_{-\ell(q)}$.  For a given $q$ we order
the basis elements for $L_{-\ell(q)}$ as follows:
$$x^{-s\mathbf d-\mathbf 1}\mathbf y^{-t\mathbf e - \mathbf 1} >
x^{-s'\mathbf d-\mathbf 1}\mathbf y^{-t'\mathbf e - \mathbf 1}$$
if and only if $s>s'$.  With respect to these ordered bases,
the matrix representing $\delta_q$ is 
$$
\begin{pmatrix}
u & v & 0 &0 & \cdots & 0 & 0 \\
0 & u & v & 0& \cdots & 0 & 0 \\
0 & 0 & u & v & \cdots & 0 & 0  \\
\vdots & \vdots & \vdots  &\ddots & \ddots & \vdots & \vdots \\
0 & 0 & 0 & 0 & \cdots & u & v \\
\end{pmatrix}_{(q+1)\times (q+2)}.
$$
By Lemma \ref{matrixlemma}, 
if $J=\cap_{q\ge 1} \ann_T \operatorname{coker}\delta_q$ then $\dim T/J=\dim
T=2$.  As $\operatorname{coker}\delta_q$ is a direct $T$-summand of $H^n_I(R/fR)$, we have
$\ann_T H^n_I(R/fR)\subseteq J$.  This completes the proof.
\qed

\begin{lem} \label{notfg} Under the assumptions of Lemma \ref{keylemma},
$\Hom_R(R/I,H^n_I(R/fR))$ is not finitely generated as an $R$-module.  Consequently,
$\Hom_R(R/I,H^n_I(R/fR))_{k}\neq 0$ for infinitely many $k$.

\end{lem}

{\it Proof:}  Suppose $\Hom_R(R/I,H^n_I(R/fR))$ is finitely generated.
By Lemma 3.5 of \cite{MV} we have that 
$I+\ann_R H^n_I(R/fR)$ is $P$-primary.  (One should note that 
the hypothesis in \cite[Lemma 3.5]{MV} that the ring
be complete  is not necessary.) This implies that $\ann_R H^n_I(R/fR)\cap T=\ann_T 
H^n_I(R/fR)$ is $m$-primary, contradicting Lemma \ref{keylemma}.
\qed

We now give the proof of our main result:
\medskip

{\it Proof of Theorem \ref{mainresult}:}  By Corollary \ref{supportmax}, $\supp H^n_I(R/fR)=
\{P\}$.  Thus, $\Hom_R(R/I,H^n_I(R/fR))_k$ has finite length as a $T$-module for all $k$ and
is nonzero for infinitely many $k$ by Lemma \ref{notfg}.  Consequently,
$$\Hom_R(R/P,H^n_I(R/fR))_{k}=\Hom_T(T/m,\Hom_R(R/I,H^n_I(R/fR))_{k})$$ is nonzero for 
infinitely many $k$.  Hence $$\starsoc_R(H^n_I(R/fR))=
\Hom_R(R/P,H^n_I(R/fR))$$ is not finitely generated.
\qed

\end{document}